\documentclass{article}

\usepackage[margin=1in]{geometry}

\title{A Statistical Framework for Domain Shape Estimation in Stokes Flows}

\author{J. Borggaard, N. Glatt-Holtz, J. Krometis  \\
  \scriptsize{emails: jborggaard@vt.edu, negh@tulane.edu, jkrometis@vt.edu} }
\date{}

\usepackage{amsfonts, amssymb, amsmath, amsthm, multicol}
\usepackage[colorlinks=true, pdfstartview=FitV, linkcolor=blue,
            citecolor=blue, urlcolor=blue]{hyperref}
\usepackage{xcolor}
\usepackage{accents}
\usepackage{comment}
\usepackage{float}
\usepackage{graphicx}
\usepackage{multirow}
\usepackage{parcolumns}
\usepackage{subcaption}
\usepackage{textcomp}
\usepackage{mathrsfs}
\usepackage{mathtools}
\usepackage{dsfont}
\usepackage{enumerate}
\usepackage[capitalize,nameinlink,noabbrev]{cleveref}
\usepackage{tabularx}

\usepackage[section]{algorithm}      
\usepackage[noend]{algpseudocode}    

\usepackage[textsize=tiny]{todonotes} 

\usepackage{gensymb} 

\usepackage{cjhebrew}

\usepackage{pdfsync}
\numberwithin{equation}{section}

\newtheorem*{Thm*}{Theorem}

\newcommand{\RR}{\mathbb{R}}
\newcommand{\x}{\mathbf{x}}
\newcommand{\y}{\mathbf{y}}

\newcommand{\E}{\mathbb{E}}

\newcommand{\DD}{\mathcal{D}}
\newcommand{\pth}{\phi}
\newcommand{\pr}{| \cdot |}
\newcommand{\bdf}{b} 

\newcommand{\ib}{\Gamma^i_{\bdf}} 
\newcommand{\ob}{\Gamma^o} 
\newcommand{\ir}{r_{\text{min}}} 
\newcommand{\mr}{r_{\text{max}}} 
\newcommand{\oR}{R} 
\newcommand{\lnk}{c} 
\newcommand{\bU}{\mathbf{u}}

\newcommand{\rtrt}{\bar{\omega}}

\newcommand{\pS}{\theta}
\newcommand{\bnh}{\hat{\mathbf{n}}}

\newcommand{\Obs}{\mathcal{O}}
\newcommand{\qS}{\mathfrak{Q}} 

\newcommand{\visc}{\nu}
\newcommand{\src}{q}
\newcommand{\data}{y}
\newcommand{\sv}{S}

\newcommand{\Vort}{\xi}

\newcommand{\fMap}{\mathcal{G}}
\newcommand{\Data}{\mathcal{Y}}
\newcommand{\noise}{\eta}
\newcommand{\mpr}{\mu_0}        
\newcommand{\mps}{\mu}          
\newcommand{\cov}{\mathcal{C}}  
\newcommand{\mcmcsamp}{\bdf}
\newcommand{\mcmccand}{\tilde{\bdf}}

\newcommand{\Exp}{\mathbb{E}} 

\begin{document}
\markboth{}{}

\maketitle

\begin{abstract}
  We develop and implement a Bayesian approach for the estimation of the shape 
  of a two dimensional annular domain enclosing a Stokes flow from sparse and noisy observations of the enclosed fluid.  
  Our setup includes the case of direct observations of the flow field as well as the measurement of concentrations of a 
  solute passively advected by and diffusing within the flow.   
  Adopting a statistical approach provides estimates of uncertainty in the shape due both to the non-invertibility of the forward map and to error in the measurements. When the shape represents a design problem of attempting to match desired target outcomes, this ``uncertainty'' can be interpreted as identifying remaining degrees of freedom available to the designer. 
  We demonstrate the viability of our framework on three concrete test problems.
  These problems illustrate the promise of our framework for applications while providing
  a collection of test cases for recently developed Markov Chain Monte Carlo (MCMC) algorithms designed 
  to resolve infinite dimensional statistical quantities.
\end{abstract}

{\noindent \small
   }

\setcounter{tocdepth}{1}
\tableofcontents

\newpage

\section{Introduction}

A variety of physical measurement problems may be formulated in terms of determining an infinite dimensional parameter 
$b$ appearing in a collection of partial differential equations from the sparse and noisy observations of the resulting
solutions.  Typically such problems are ill-posed in that it is impossible to exactly recover the unknown $b$ from the available data.
For such problems a suitable Bayesian approach provides a principled and comprehensive probabilistic framework for the estimation of $b$, \cite{stuart2010inverse}.  Such a Bayesian methodology is viable due to recent developments in computational statistics \cite{cotter2013mcmc}
and advances in high performance computing.
This work is devoted to demonstrating the viability of a Bayesian framework for the estimation of the shape of a domain enclosing a steady Stokes fluid flow with a passive scalar.

Although designing shapes of objects in fluids has a rich history, their treatment within a mathematical framework can be traced to the work of Pironneau~\cite{pironneau1973OptimumProfilesStokes} that studied optimal shapes of objects in Stokes flows.  The traditional framework is to set up a constrained optimization problem to minimize an objective function that depends on the shape of the domain through the solution to partial differential equations that model the fluid.  There have been a number of notable works that discuss computation of required derivatives of shape parameters~\cite{haslinger1987SensitivityAnalysisOptimal,sokolowski1992introduction,borggaard1993SensitivityCalculations2D,jameson1995OptimumAerodynamicDesign,delfour2011shapes,hiptmair2015ComparisonApproximateShape} and different optimization approaches, e.g.~\cite{kelley1991MeshIndependenceNewtonlike,arian1997CouplingAerodynamicStructural,arian1995ShapeOptimizationOneshot,borggaard1997PDESensitivityEquation,alexandrov2000OptimizationVariablefidelityModels,hinze2009OptimizationPDEConstraints}.  Other significant rigorous literature related to shape optimization include \cite{sokolowski1992introduction}, \cite{delfour2011shapes}, and \cite{schulz2014riemannian}.

As is typical of many high-dimensional inverse problems, the shape of the domain is unlikely to be specified uniquely by the lower-dimensional data. In fact, there may be many classes of shapes that yield similar matches to the target criteria. 
From the perspective of an inverse problem -- i.e., where the shape is unknown and the practitioner is trying to infer it from the data -- this uncertainty provides key insight into how much information the data provides.
The shape optimization problem can also be viewed from the perspective of design, where the data represents outcomes that are desired from the system -- such as the fluid behavior that we describe in this paper -- and a shape is sought to meet those outcomes. In this case, the ``uncertainty'' represents degrees of freedom that may still be available to the designer that can be used to meet other objectives (e.g., aesthetics, secondary technical constraints, etc). 
All of these considerations make it desirable to treat the shape estimation problem in a manner such that the results tell us what the data does and does not dictate about the unknown shape.

The Bayesian approach to inverse problems \cite{kaipio2006statistical,dashti2017bayesian,stuart2010inverse,gelman2014bayesian} addresses this need by treating both the unknown parameter and the target data as statistical quantities. The distribution of the former (known as the \emph{prior}) might be based on historical data or theoretical knowledge about the parameter structure, while the latter (known as the \emph{likelihood}) is typically based on accuracy of the sensors that measured the data (in the ``uncertain data'' case) or the allowed tolerances on target values (in the ``design'' case). The output of the framework is then a probability distribution $\mu$ on the parameter space known as the \emph{posterior} that tends to have largest probability where the parameter matches both the prior and the likelihood well, or one especially well, and lowest probability where the parameter has a significant mismatch with one or the other.

To illustrate a statistical approach to shape estimation, in this paper we consider the problem of estimating the unknown inner boundary of a cylindrical mixer to achieve desired behavior of either the fluid inside the mixer or a substance moving within the fluid. The fluid movement is modeled via the steady two-dimensional Stokes equations, a low-Reynolds number approximation to the Navier-Stokes equations. The substance in the fluid is considered to be a passive scalar -- i.e., its behavior is a byproduct of the fluid movement but does not in turn affect the fluid velocity or pressure -- and modeled by coupling the Stokes flow to the advection-diffusion equation, which models the scalar concentration in the presence of convective and diffusive forces. We adopt the Bayesian methodology described above to estimate uncertainty in the inner boundary and present the results for three problems with different target data, one based on the flow itself and the other two based on the scalar.   Regarding the prior we adopt a Gaussian methodology which calibrates the degree of smoothness/roughness of inner boundary in terms of a covariance operator involving powers of the inverse of Laplacian operator on the circle.

It is worth emphasizing that effectively resolving posterior $\mu$ resulting from the statistical inversion procedure used in our problem essentially relies on recent algorithmic advances in computational statistics \cite{beskos2008mcmc, cotter2013mcmc, Beskosetal2011, glatt-holtz2020acceptreject,glattholtz2022on}.  Indeed for the three test problems we consider here we demonstrate the efficacy of the so-called preconditioned Crank-Nicolson (pCN) algorithm for the accurate resolution of a variety of observables for $\mu$.   However there is considerable scope for the design and optimal tuning of algorithmic parameters for infinite-dimensional sampling algorithms as we illuminated in our recent works \cite{glatt-holtz2020acceptreject,glattholtz2022on}. Thus these test problems serve as significant benchmarks for testing the relative performance and the optimal tuning of parameters in a variety of sampling methods.   In particular in \cite{glattholtz2022on} the test case \cref{sec:svsector} serves as an initial demonstration of a novel `multiproposal' method that is tailored for parallel computing architectures in widespread use.

To our knowledge, little work has been done to date on the application of the Bayesian framework to shape optimization problems. 
Two recent exceptions are \cite{iglesias2014wellposed} and \cite{iglesias2016bayesian}, which develop a Bayesian level set method and related numerical algorithms for PDE-based problems with unknown domain boundaries.
Also, \cite{kawakami2020StabilitiesShapeIdentification} provides theoretical results for the stability of the statistical approach to estimating the shape of a cavity in a heat conductor from the evolution of the heat inside as modeled by the heat equation. 
Other works considering shape estimation from a statistical (but not Bayesian) perspective include~\cite{huyse2001AerodynamicShapeOptimization,koumoutsakos2006FlowOptimizationUsing}. 
Works that apply a traditional/deterministic approach to fluids problems similar to the one that we consider here include \cite{eggl2020shape}, which seeks the boundary shape to maximize mixing of a binary fluid in a system governed by the Navier-Stokes equations, \cite{cohen2018ShapeHolomorphyStationary}, which proves shape holomorphy for the map from boundary parameters to PDE solution for Stokes and Navier-Stokes flows, and \cite{hu2021BoundaryControlOptimal}, which proves the existence of the optimal boundary control to mix a similar Stokes-passive scalar system. 

The remainder of this paper is organized as follows: In \cref{sec:formal:setup}, we describe how the boundary of the domain is parameterized, the PDEs that we solve on the resulting domain, the observations of those solutions that we consider, and how the shape optimization problem is formulated as a statistical inverse problem. \cref{sec:numerical:method} provides the details of the numerical methods that we use to map the parameter to the boundary, to solve the PDEs, and to sample from the posterior measure constituting the answer to the inverse problem. In \cref{sec:example:problems}, we outline application of the statistical approach to three example problems with different features. Finally, some concluding thoughts and next steps are provided in \cref{sec:conclusions}.

\section{Formal Set-Up}
\label{sec:formal:setup}

In this section we lay out the formal set up that allows us to
formulate the domain shape estimation problem in a Bayesian
statistical inversion framework,  \cite{kaipio2006statistical,dashti2017bayesian}.

\subsection{The Boundary to Stokes Flow Map}
\label{sec:forward:map}

\begin{figure}[ht!]
	\centering
	\includegraphics[width=0.6\linewidth]{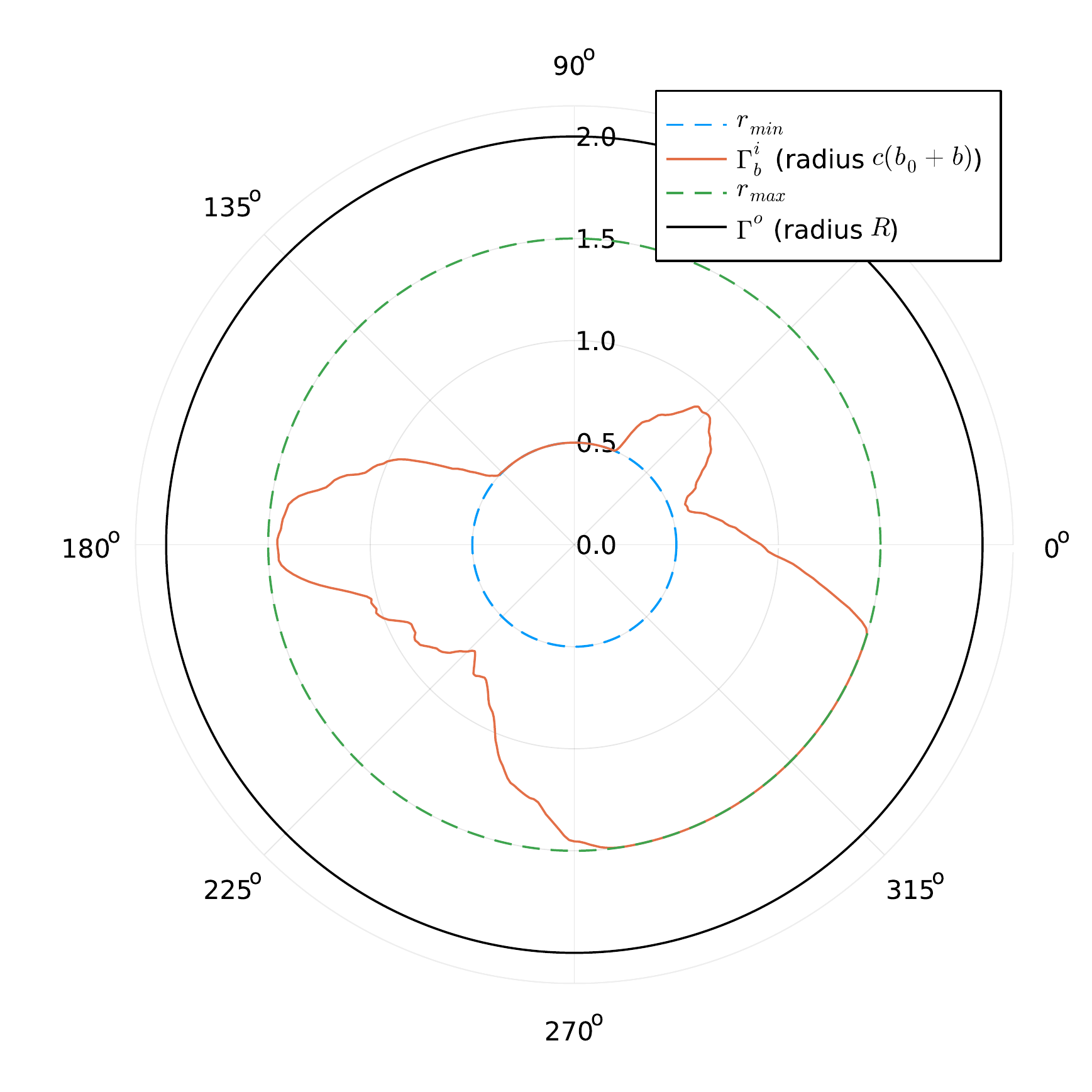}
	\caption{The inner and outer boundaries $\ib,\ob$ and radius constraints $\ir$ and $\mr$.}\label{fig:radii}
\end{figure}

We consider a class of Stokes flows on subsets of the domain
\begin{align}
  \DD := \{\x \in \RR^2 | \ir < |\x| \leq \oR\},
  \label{eq:com:domain}
\end{align}
defined so the flow always contains the subdomain
\begin{align}
  \DD^0 := \{\x \in \RR^2 | \mr \leq |\x| \leq \oR\},
  \label{eq:com:sub:domain}
\end{align}
for some $0< \ir < \mr < \oR$.  We denote by
$(\pth, \pr) : \DD \to [0,2\pi) \times [\ir,\oR]$ the conversion to polar
coordinates.  We next fix a `clamping' function $\lnk: \RR \to (\ir, \mr)$ to be a
smooth strictly increasing function with
\begin{align*}
  \lim_{t \to -\infty} \lnk(t) = \ir,
  \quad
  \lim_{t \to \infty}\lnk(t) = \mr.
\end{align*}
Many forms of $\lnk$ may be used; our choice is described in detail in \cref{sec:clamp}, see \eqref{eq:clamp:poly}. 
Given any $\bdf:\RR \to \RR$ which is sufficiently smooth and 
$2\pi$ periodic we consider the domain
\begin{align}
  \DD_{\bdf} := \{\x \in \RR^2 |  \lnk(\bdf_0+\bdf(\pth(\x))) \leq |\x| \leq \oR\},
  \label{eq:b:dom}
\end{align}
where $\bdf_0$ is a fixed mean radius and denote the inner and outer boundaries of $\DD_{\bdf}$, respectfully, as
\begin{align}
  \ib := \{\x \in \RR^2 | \lnk(\bdf_0+\bdf(\pth(\x))) = |\x| \},
  \quad
  \ob := \{ \x \in \RR^2 | |\x| = \oR\}.
\end{align}
The relationship between the inner boundary and various radii is visualized in \cref{fig:radii}.

For each $\bdf:\RR \to \RR$ sufficiently smooth and periodic we
consider the following steady Stokes problem:
\begin{align}
  \visc \Delta \bU = \nabla p, \quad \nabla \cdot \bU = 0
  \quad \text{ in } \DD_\bdf,
  \label{eq:Sk:bulk}
\end{align}
where $\bU$ is the velocity and $p$ is the pressure and $\visc > 0$ is
the kinematic viscosity, a known parameter in our problem.  For Stokes flows, $\visc$ does not affect the velocity field, but does scale the pressure and fluid stress terms that could be measured.
Equation \eqref{eq:Sk:bulk} is supplemented with the boundary conditions
\begin{align}
  \bU = \rtrt \x^\perp \text{ on } \ob
  \quad  \text{ and } \quad 
  \bU = 0 \text{ on } \ib,
  \label{eq:Sk:bnd}
\end{align}
where $\rtrt \in \RR$ is another known scalar parameter representing a
rotation rate on the outer boundary and $\x^\perp = (x_2, -x_1)$.

In several of the examples below, we couple
\eqref{eq:Sk:bulk}--\eqref{eq:Sk:bnd} to a steady advection
diffusion equation
\begin{align}
  \bU \cdot\nabla \pS =  \kappa \Delta \pS + \src
  \quad \text{ in } \DD_\bdf
  \label{eq:AD:steady}
\end{align}
where $\kappa > 0$ is the diffusivity parameter for the passive solute
$\pS$ and $\src$ is a volumetric source which we typically assume to be
compactly supported in the common subdomain $\DD^0$ as in
\eqref{eq:com:sub:domain}.  Regarding boundary conditions for
\eqref{eq:AD:steady} we consider a homogeneous Dirichlet outer
condition and an insulated inner boundary. Namely we suppose that
\begin{align}
  \pS = 0  \text{ on } \ob
  \quad  \text{ and } \quad 
  \nabla  \pS \cdot \bnh = 0 \text{ on } \ib,
  \label{eq:AD:steady:BC}
\end{align}
where $\bnh$ is the outward pointing normal to the inner boundary
$\Gamma^i_\bdf$.

\subsection{Observational Procedures}
\label{sec:obs:procedure}

We will consider the following observational procedures
$\Obs = \Obs(\bU, \pS, \DD_b)$ for solutions of
\eqref{eq:Sk:bulk}--\eqref{eq:Sk:bnd}, sometimes involving
\eqref{eq:AD:steady}--\eqref{eq:AD:steady:BC}. These are just examples that might be computed from the solutions $\bU$ and $\pS$; many other quantities of interest could of course be used for $\Obs$. 

\subsubsection*{Case 1: Volumetric averages of the Vorticity}

The first observation procedure involves spatially averaged
measurements of vorticity $\Vort := \nabla^\perp \cdot \bU$ of
solutions $\bU$ of \eqref{eq:Sk:bulk}--\eqref{eq:Sk:bnd}.  In this
case we fix $\x_1, \ldots, \x_n$ in the subdomain $\DD^0$, \eqref{eq:com:sub:domain},
common to every parameter value $\bdf$ and radii
$r_1, \ldots, r_n$ so that $B(\x_j, r_j) \subset \DD^0$
for every $j = 1, \ldots, n$.  Here $B(\x, r)$ is the ball of radius
$r$ with center $\x$.  In this case
\begin{align}
  \Obs(\bU, \pS, \DD_b) := \left(\frac{1}{\pi r_1^2} \int_{B(\x_1, r_1)}
  \nabla^\perp \cdot \bU(\x) d\x,\quad  \ldots,\quad   \frac{1}{\pi r_n^2}
  \int_{B(\x_n, r_n)}\nabla^\perp \cdot \bU(\x) d\x \right).
  \label{eq:Vol:ave:msr}
\end{align}

\subsubsection*{Case 2: Domain average of the Scalar Variance}

In this second situation we measure the global average of the scalar
variance over the domain $\DD_\bdf$. Scalar variance is a key measure of mixing given by
\begin{equation}\label{eq:sv}
    \sv = \int_{\DD_\bdf} \left(\pS(\x) - \bar{\pS}\right)^2 \,d\x
    \quad \text{ where } \bar{\pS} = \frac{1}{|\DD_\bdf|}\int_{\DD_\bdf} \pS(\y) d \y
\end{equation}
Here $|\DD_\bdf| = \int_{\DD_\bdf} d\x$ is the total area of
$\DD_\bdf$
so that $\bar{\pS}$ is the mean value of $\pS$ over the domain $\DD_\bdf$. Low scalar variance indicates that the scalar has, via advection and diffusion, been spread evenly throughout the domain (e.g., $\sv=0$ for $\pS$ constant in $\x$), while higher scalar variance suggests that the scalar has been ``trapped'' in a subregion. To avoid biasing the choice of shape toward domains with uniformly large or small radii, in this case we normalize by volume. That is, we make a single
measurement of the solution \eqref{eq:Sk:bulk}--\eqref{eq:AD:steady:BC} given by
\begin{align}
  \Obs(\bU, \pS, \DD_\bdf)
  := \frac{1}{|\DD_\bdf|} \int_{\DD_\bdf} \left( \pS(\x) - \bar{\pS} \right)^2 d\x,
  \label{eq:scalar:var:msr}
\end{align}
where  $\bar{\pS}$ is the mean value of $\pS$ as defined above in \eqref{eq:sv}. 

\subsubsection*{Case 3: Sectorial values of the average Scalar Variance}

A third observation procedure involves dividing $\DD_\bdf$ into angular
sectors as
\begin{align}\label{eq:sectors}
  \qS_j(\DD_\bdf) = \{ \x \in \DD_\bdf | \pth_{j -1} < \pth(\x)  \leq \pth_{j} \}
\end{align}
where $0 = \pth_0 < \pth_1 < \cdots < \pth_{n-1} < \pth_{n} = 2\pi$.  Then
we define
\begin{align}
  \Obs(\bU, \pS, \DD_\bdf) := \left(
  \frac{1}{|\DD_\bdf|} \int\limits_{\qS_1(\DD_\bdf)} \left( \pS(\x) - \bar{\pS} \right)^2 d\x,\quad
  \ldots,   \quad
  \frac{1}{|\DD_\bdf|}\int\limits_{\qS_n(\DD_\bdf)} \left( \pS(\x) - \bar{\pS} \right)^2 d\x 
  \right).
 \label{eq:scalar:sect:var:msr}
\end{align}
where, again, $|\DD_\bdf| = \int_{\DD_\bdf} d\x$ is the total area of
$\DD_\bdf$ and $\bar{\pS} = \frac{1}{|\DD_\bdf|}\int_{\DD_\bdf} \pS(\y) d \y$ is the mean value of $\pS$ over $\DD_\bdf$. 
Note that the observations are normalized by the total area of the domain so that they sum to the domain average scalar 
variance \eqref{eq:scalar:var:msr}. The domains are therefore not incentivized to be large or small as a whole, but radii 
may tend to be large or small within a given sector.

\subsection{Formulation as a statistical inverse problem}
\label{sec:stat:inv:background}

As we sketched in the introduction our goal is develop a principled
approach for estimating the shape of the domain $\DD_\bdf$ from sparse
observations of the associated solutions of \eqref{eq:Sk:bulk} and/or \eqref{eq:AD:steady}.  In other words
we would like to estimate the parameter $\bdf$ in \eqref{eq:b:dom} 
from the sort of measurements described
in \cref{sec:obs:procedure}.
For this purpose we leverage a Bayesian statistical inversion
framework introduced in e.g. \cite{kaipio2006statistical} and
developed in the infinite dimensional, PDE-constrained setting in
\cite{stuart2010inverse, dashti2017bayesian}.  

We recall this Bayesian formalism here as follows.  
Suppose are trying estimate an unknown `parameter'
$\bdf$ which we assume is an sitting in a separable Hilbert space $H$.  Each value of this
parameter produces a finite number of discrete measurements through a
forward map $\fMap: H \to \RR^n$ which typically involves the solution
of a partial differential equation and which depends on $\bdf$ in
a complicated, possibly nonlinear fashion.  While, for a given observation $\Data$, we would like
to, in effect, `invert' $\fMap$ by setting $\bdf = \fMap^{-1}(\Data)$ such a direct inversion is
often impossible in  principle and anyway is at least
intractable in practice.  Complicating matters further the measurements in $\Data$ 
may be subject to a degree of observational uncertainty which needs to
be accounted for as well.

A Bayesian solution to this problem is consider the measurement
model
\begin{align}
  \Data = \fMap(\bdf) + \noise
  \label{eq:stat:inverse:model}
\end{align}
where $\noise$ is an (additive) observational noise which we suppose is a random
quantity with a known distribution $\gamma \in Pr(\RR^n)$.
Information on $\bdf$ garnered from the observation of $\Data$ is
supplemented with other external information on the unknown parameter
$\bdf$ given in the form of a prior probability distribution
$\mpr \in Pr(H)$.  Bayes' theorem then yields a posterior probability
distribution $\mps$ which provides a comprehensive model of our
uncertainties in $\bdf$.  Indeed assuming a Gaussian observational error, namely that $\noise \sim \gamma = N(0, \Sigma)$ for a symmetric positive definite covariance $\Sigma$, in \eqref{eq:stat:inverse:model}
we obtain a posterior measure of the `explicit' form
\begin{align}
  \mps(d \bdf) =
  \frac{1}{Z} \exp\bigl( - \frac{1}{2} | \Sigma^{-1/2}(\Data
  - \fMap(\bdf))|^2\bigr) \mpr(d \bdf).
  \label{eq:post:form}
\end{align}
See e.g. \cite{dashti2017bayesian},
\cite{borggaard2020bayesian} for a formulation of Bayes' theorem needed
in the generality considered here.

Gaussian measures provide a typical choice for the prior $\mpr$ particularly 
when the unknown parameter is infinite dimensional.  This is
a natural way to, for example, specify the degree of underlying spatial
regularity of the infinite dimensional parameter $\bdf$ as in
\eqref{eq:b:dom}; in our current context this corresponds to specifying the degree of
roughness of the inner boundary of $\DD_\bdf$.  
Note this Gaussian setting is also desirable given the extensive theory for
Gaussian measures of functional spaces, cf. \cite{DPZ2014}.  Crucially on 
the basis of this theoretical foundation recent developments for Markov Chain Monte Carlo (MCMC) sampling 
\cite{beskos2008mcmc, cotter2013mcmc, Beskosetal2011, glattholtz2022on} provide sampling
methods that partially beat the curse of dimensionality in the setting of Gaussian priors.  
Indeed it is a Gaussian prior formulation that
allows the meaningful numerical resolution of statistics from
statistical posteriors in the examples we consider below in
\cref{sec:example:problems}.

Let us briefly recall some elements of the infinite dimensional Gaussian 
theory as suits our purposes here.   Given a separable Hilbert space
$H$, a Borel probability measure $\nu \in Pr(H)$ is said to be Gaussian
if $\ell^* \nu$ is a (one dimensional) normal distribution for any bounded linear
functional $\ell$.  Here we recall that $\ell^* \nu(A) = \nu(\ell^{-1}(A))$ is 
the pushforward of $\nu$ under $\ell$, that is to say if $X \sim \nu$ then
$\ell(X) \sim \ell^*\nu$.  Now for any such Gaussian measure $\nu$ there 
is a uniquely defined `mean' $m \in H$ and covariance $\cov$, a symmetric, positive
and trace class operator on $H$ such when $X \sim \nu$
\begin{align*}
 \E \langle X, f \rangle =  \langle m, f \rangle, 
 \quad \E (\langle X-m, f \rangle \langle X-m, g\rangle)  =  \langle \cov f, g \rangle, 
\end{align*}
for any $f, g \in H$.  In general for any bounded $\cov$ that is symmetric, positive
and trace class we may obtain an orthonormal basis $\{e_k\}_{k \geq 1}$ of $H$ which are 
eigenfunctions of $\cov$, $\cov e_k = \lambda_k e_k$.  For any such eigenbasis $\{e_k\}_{k \geq 1}$
it is not hard to see that by defining
\begin{align}
	X =m +  \sum_{k = 1}^\infty \sqrt{ \lambda_k} e_k \xi_k \quad \text{ for } 
	\xi_k \text{ taken as i.i.d } N(0,1),
	\label{eq:KL}
\end{align}
we obtain a random variable distributed as $\nu$.  The representation \eqref{eq:KL}
is often referred to as a Karhunen-Lo\'eve expansion of $\nu$ and provides a means of simulating
such distributions $\nu$.

To place our inverse problem in the general setting of
\eqref{eq:stat:inverse:model}--\eqref{eq:post:form} we consider our
unknown parameter space as a Sobolev space which can be defined, for $s \geq 0$, 
as
\begin{align}
  H^s := \left\{ \bdf: [0,2\pi] \to \RR :
  \bdf(x) := 
  \sum_{k=1}^\infty (\bdf_{2k-1} \cos (k x) + \bdf_{2k} \sin
  (k x)),
  \|\bdf\|_{H^s} < \infty \right\}
  \label{eq:sob:sp}
\end{align}
where
\begin{align}
  \| \bdf\|_{H^s}^2 := 
  \frac{1}{2}\sum_{k=1}^{\infty} k^{2s} (\bdf_{2k-1}^2+\bdf_{2k}^2).
    \label{eq:sob:sp:norm}
\end{align}
In other words $H^s$ is morally the set of $2\pi$ periodic, mean zero functions with
$s$ derivatives in $L^2$.
Note the $H = H^0$ corresponds to $L^2$ and that standard Sobolev embedding results $H^s \subset C^{\tilde{s}}$ hold
whenever $s > 1/2 + \tilde{s}$, $\tilde{s} \geq 0$ (see e.g. \cite{robinson2001infinite}).
To parameterize the inner boundary, we let the components $\{\bdf_k\}_{k=1}^\infty$ be our unknown parameters. 

The forward map $\fMap(\bdf)$ for $\bdf \in H^s$ with $s > 1/2$ is specified by defining $\DD_\bdf$ according to
\eqref{eq:b:dom} and then solving for $\bU_\bdf$ according to
\eqref{eq:Sk:bulk}--\eqref{eq:Sk:bnd} and then (if necessary for the observations) subsequently
solving for $\pS_\bdf$ via \eqref{eq:AD:steady}--\eqref{eq:AD:steady:BC}.  With these
elements in hand we can then define $\fMap(\bdf) := \Obs(\bU_\bdf,
\pS_\bdf, \DD_\bdf)$ for a given observation operator $\Obs$ as described in \cref{sec:obs:procedure}. 
The numerical implementation of these steps are described in \cref{sec:forward:map}.

Regarding the prior (Gaussian) measure on these parameters, we choose
$\mpr = N(0,\cov)$ where $\cov$ is a diagonal operator such that components
$\bdf_k, k\ge 1$ are mutually independent and distributed
as 
\begin{align}
	\bdf_{2k-1},\bdf_k \sim N(0,k^{-2s-1});
	\label{eq:prior:cov:sp}
\end{align}
see \eqref{eq:KL}. Using \eqref{eq:sob:sp:norm},
this ensures that draws from the prior are almost surely elements of
$H^{s'}$ for any $s' < s$:
\begin{align}\label{eq:stokes:prior}
    \Exp \| \bdf\|_{H^{s'}}^2 
    = \Exp \left[ 
    \frac{1}{2}\sum_{k=1}^\infty k^{2s'} (\bdf_{2k-1}^2+\bdf_{2k}^2) \right]
    = 
    \sum_{k=1}^\infty k^{-1 - 2(s-s')} < \infty.
\end{align}

\section{Numerical Methods}
\label{sec:numerical:method}

Computing meaningful quantities from the posterior measure \eqref{eq:post:form} typically involves using numerical methods to sample from it. For this, we use Markov Chain Monte Carlo (MCMC) methods developed in recent years for sampling from function spaces; these methods are described in detail in \cref{sec:mcmc}. Generating the samples involves computation of the forward map $\fMap$ -- i.e., solving the Stokes \eqref{eq:Sk:bulk}--\eqref{eq:Sk:bnd} and advection-diffusion equations \eqref{eq:AD:steady}--\eqref{eq:AD:steady:BC} and computing the observations associated with these solutions; our implementation of the forward map is described in \cref{sec:forward:map}. The full solver and MCMC routines were written in the Julia numerical computing language \cite{bezanson2017julia} and are publicly available at \url{https://github.com/jborggaard/BayesianShape}.

\subsection{Monte Carlo Sampling Algorithms}\label{sec:mcmc}

As noted above, Markov Chain Monte Carlo (MCMC) provides a common method for sampling from the posterior measure \eqref{eq:post:form}. In particular, we will employ the Metropolis-Hastings algorithm, which dates from the mid-twentieth century \cite{metropolis1953equation,hastings1970monte} and uses an accept-reject mechanism to convert samples from distributions that are easy to draw pseudo-random samples from into samples from an arbitrary target measure. Much effort has been expended in recent years to adapt this method to general state spaces and in particular to cases where the unknown is a function as we describe in \cref{sec:stat:inv:background}. These include theoretical developments, casting the method in general terms and state spaces \cite{tierney1998note,neklyudov2020involutive,andrieu2020general,glatt-holtz2020acceptreject}, and algorithmic advances to develop methods well-suited to states spaces of arbitrary dimension (see, e.g., \cite{Beskosetal2011,cotter2013mcmc, glattholtz2022on}), improving practical performance. 

While high dimensional gradient-based methods such as $\infty$-Hamiltonian Monte Carlo (HMC) \cite{Beskosetal2011} or $\infty$-Metropolis-adjusted Langevin (MALA) \cite{cotter2013mcmc} (see also \cite{stuart2004conditional}) methods could be employed here in principle,
the implementation costs to employ a suitable adjoint method are significant. The reliance on a third-party, black box meshing algorithm (see \cref{sec:fem}) also complicates the use of automatic differentiation \cite{borggaard2000efficient,fournier2012ad} to compute gradients automatically. We therefore use the preconditioned Crank-Nicolson (pCN) algorithm from \cite{beskos2008mcmc, cotter2013mcmc} in the numerical studies in this paper. pCN is a gradient-free method that adapts the classical ``random walk'' Metropolis-Hastings method to the infinite-dimensional setting with Gaussian prior $\mpr = N(0,\cov)$. The algorithm is summarized in \cref{alg:pcn}.   

Note that, given the difficulties of computing gradients for this problem several recently identified alternatives MCMC methods present themselves.  In \cite{glatt-holtz2020acceptreject}
a class of $\infty$HMC algorithms are derived which provide a bias free methodology using essentially any reasonable approximation of the gradient.  
On the other hand in  \cite{glattholtz2022on}  we recently derived a `multi-proposal' variation of pCN
in which holds promise as a different `gradient-free' alternative to $\infty$HMC and which is well suited for parallel computational architectures.  In \cite{glattholtz2022on}, we present some preliminary convergence results for this new mpCN algorithm using the problem described in \cref{sec:svsector} as one of our test cases.

In the numerical examples in \cref{sec:example:problems}, we adaptively tune the ``step size'' parameter $\rho$ to target an acceptance rate of roughly $23\%$ which has been suggested as optimal for some problems \cite{roberts2001optimal}, with the adaptation window starting as $[0.5,2]$ and vanishing linearly to zero as the run progresses.  See, e.g., \cite{roberts2009examples} or \cite[Section 12.2]{gelman2014bayesian} for discussions of adaptive MCMC parameter tuning.

\begin{algorithm}
\caption{Preconditioned Crank-Nicolson (pCN).}\label{alg:pcn}
\begin{algorithmic}[1]
\item Given free parameter $\rho \in [0,1)$ and initial sample $\mcmcsamp^{(j)}$
\item Propose $\mcmccand = \rho \mcmcsamp^{(j)} + \sqrt{1-\rho^2} \xi^{(j)}$, $\xi^{(j)} \sim \mpr = N(0,\cov)$ (see \eqref{eq:KL}, \eqref{eq:prior:cov:sp}).
\item Set $\mcmcsamp^{(j+1)} = \mcmccand$ with probability
  $\min\left\{1,\exp\left(\Phi\left(\mcmcsamp^{(j)}\right)
      - \Phi\left(\mcmccand\right) \right)\right\}$, otherwise $\mcmcsamp^{(j+1)} = \mcmcsamp^{(j)}$
\end{algorithmic}
\end{algorithm}

\subsection{Computing the Forward Map $\fMap$} \label{sec:forward:map}

To compute the forward map $\fMap$ numerically from a given (finite) array of components $\{\bdf_k\}_{k=1}^K$ which parameterize the inner boundary shape $b$, we need to solve the Stokes \eqref{eq:Sk:bulk}--\eqref{eq:Sk:bnd} and advection-diffusion \eqref{eq:AD:steady}--\eqref{eq:AD:steady:BC} equations. For this, we use a finite element solver, which requires a mesh of the domain $\DD_b$, for which we use Gmsh \cite{geuzaine2009gmsh}, a popular meshing software. This in turn required representing the boundaries in a format that could be passed to Gmsh; for this, we convert the Fourier representation to B-splines. Our procedure for computing the forward map therefore involves the following steps:
\begin{enumerate}
\item Compute the Fourier expansion \eqref{eq:sob:sp}, 
\item Compute the optimal B-spline approximation to the Fourier expansion (Section~\ref{sec:bspline}),
\item Clamp via $\lnk$ (defined below in \eqref{eq:clamp:poly}) to get a B-spline representation of the inner boundary $\ib$ (Section~\ref{sec:clamp}),
\item Mesh the interior of the domain $\DD_\bdf$ with Gmsh \cite{geuzaine2009gmsh} (Section~\ref{sec:fem}), 
\item Solve the
Stokes PDE \eqref{eq:Sk:bulk}--\eqref{eq:Sk:bnd}, and, if required for the observations, the advection-diffusion PDE \eqref{eq:AD:steady}--\eqref{eq:AD:steady:BC}, via a finite element solver, and finally,
\item Compute the observations as described in \cref{sec:obs:procedure}.
\end{enumerate}

These steps are described, as necessary, in the following subsections.

\subsubsection{Fourier Expansion and B-Spline Approximation} \label{sec:bspline}
In this subsection, we describe how we compute the optimal B-spline approximation to the inner radius from a given Fourier expansion $r_F(x;\bdf^{(j)})$ generated from the parameter array $\bdf^{(j)}$.  To do this, we partition the periodic domain $[0,2\pi)$ into $n_B$-uniformly spaced intervals $[x_i,x_{i+1}], i=1,\ldots,n_B$ with $x_i=x_{n_B+i}$ that generate a basis of $n_B$ cubic B-splines, cf.~\cite{evans2009n-widths},
\begin{align}
\label{eq:bSpl:basis}
  \mbox{span}\{ B_{1,3}, B_{2,3}, \ldots, B_{n_B,3} \}
\end{align}
defined recursively using the relationships
\begin{displaymath}
  B_{i,0}(x) = \left\{\begin{array}{cl}
  1 & \mbox{if }x_i\leq x\leq x_{i+1} \\
  0 & \mbox{otherwise}
  \end{array} \right. \quad \mbox{and} \quad
  B_{i,p}(x) = \frac{x-x_i}{x_{i+p}-x_i}B_{i,p-1}(x) + \frac{x_{i+p+1}-x}{x_{i+p+1}-x_{i+1}}B_{i+1,p-1}(x)
\end{displaymath}
for $p=1,2,\ldots$.

The coefficients of the B-spline expansion $r_B(x) = \sum_{i=1}^{n_B} a_i B_{i,3}(x)$ are found as the best $L_2$-approximation to $r_F$ using the projection theorem:
\begin{displaymath}
  \left\langle r_F-r_B, B_{i,3} \right\rangle_{L_2} = 0, \qquad i=1, \ldots, n_B.
\end{displaymath}
Due to the local support of the cubic B-splines, this leads to a symmetric Toeplitz system with 7 diagonal entries and a right-hand-side determined by combining the coefficients $\bdf_k$ with integrals of products of polynomials and trigonometric functions that can be analytically precomputed.  The result is a set of B-Spline coefficients $a_i$ along with an error in the approximation using the $n_B$ intervals. The number of B-splines required to accurately represent $r_F$ by $r_B$ depends on the number of Fourier components used in $r_F$ and the regularity of the boundary imposed by the prior.  For the numerical examples in \cref{sec:example:problems}, the $L^2$ error in this B-spline approximation ranged from roughly 0.005 to 0.02.

\subsubsection{The ``Clamping'' Function $\lnk$} \label{sec:clamp}

\begin{figure}[h]
    \centering
    \includegraphics[width=0.48\textwidth]{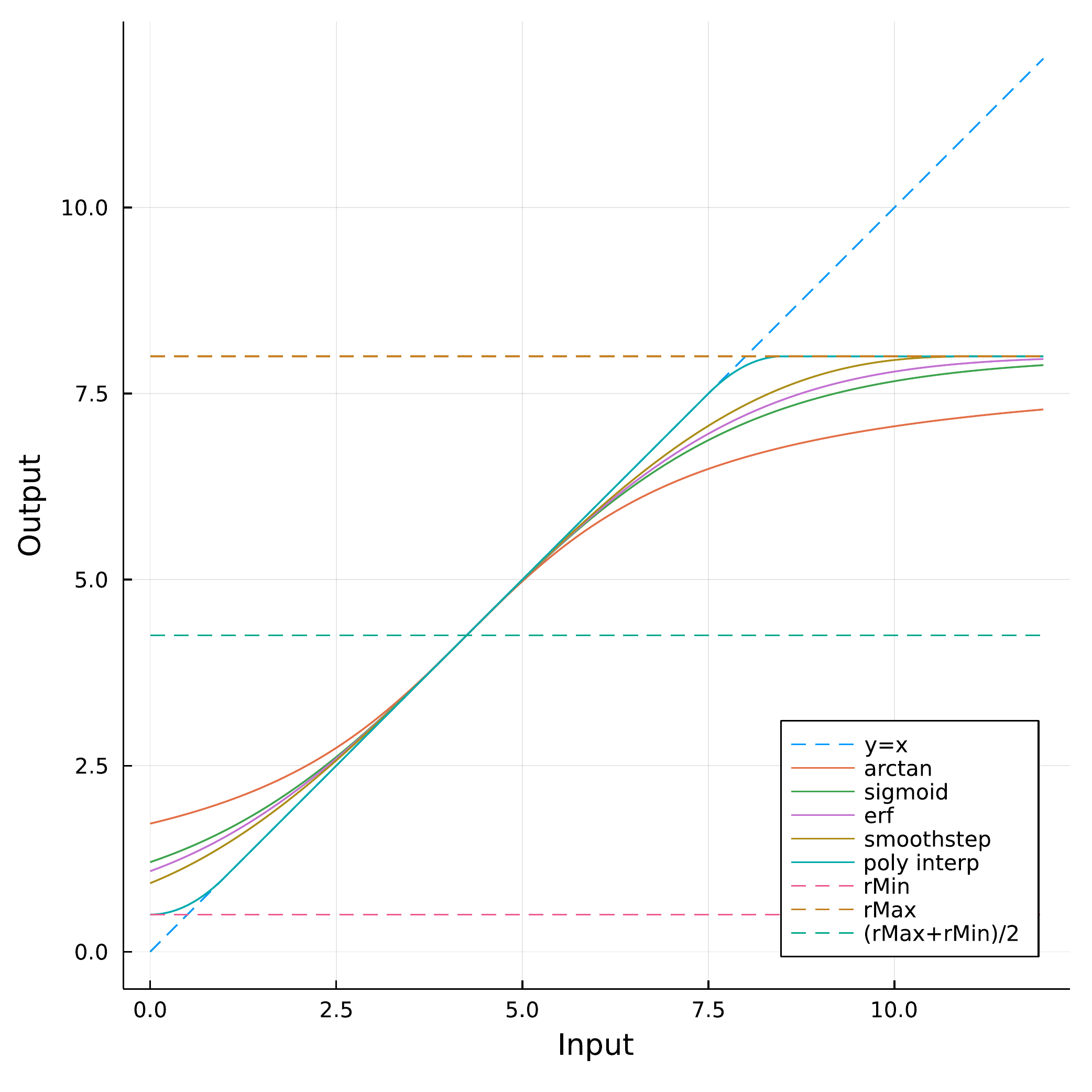}
    \hfill
    \includegraphics[width=0.48\textwidth]{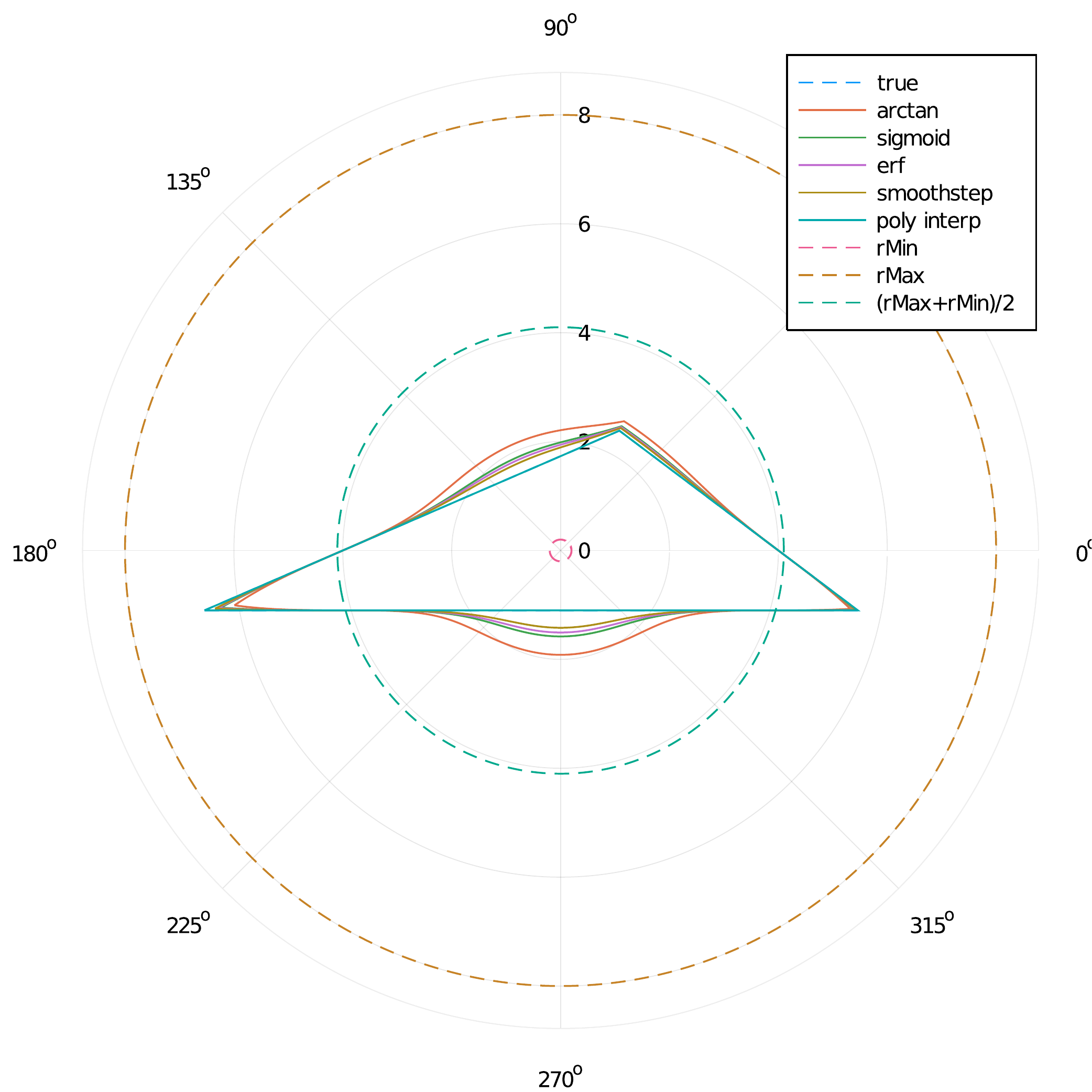}
    \caption{Example clamping functions $\lnk$. Left: The effect of $\lnk$ on different radius values; the polynomial interpolant \eqref{eq:clamp:poly} was chosen for numerical experiments because it provides a smooth function with sufficient degrees of freedom to limit radius changes to close to the boundary. Right: Various $\lnk$ applied to a triangular inner boundary; the polynomial interpolant lies directly on top of the ``true'' (unmanipulated) boundary.}
    \label{fig:link}
\end{figure}

There are many standard choices in computing literature for ``clamping'' functions -- smooth, non-decreasing functions used to map an input from $(-\infty,\infty)$ to a bounded domain $[\ir,\mr]$. Examples include sigmoid, arc tangent, the error function $\text{erf}$, and the ``smoothstep'' function. However, experimentation showed that the impact of these functions was too global for this application, i.e., these functions began smoothing the boundary too far away from the constraints $\ir$ and $\mr$. We therefore developed a custom polynomial interpolant that adjusted the radius only within a given radius $\epsilon$ of either boundary, while providing smooth transitions at the interfaces:
\begin{align}\label{eq:clamp:poly}
    \lnk(t) = \begin{cases}
      \ir & t \le \ir - \epsilon \\
      \ir + \frac{1}{4\epsilon}\left(t-\ir+\epsilon\right)^2 & \ir-\epsilon < t < \ir+\epsilon \\ 
      t & \ir+\epsilon \le t \le \mr-\epsilon \\ 
      \mr-\frac{1}{4\epsilon}\left(t-\mr-\epsilon\right)^2  & \mr-\epsilon < t < \mr+\epsilon \\ 
      \mr & \mr + \epsilon \le t
    \end{cases}.
\end{align}
\cref{fig:link} compares the custom interpolant with the more standard choices of clamping functions. The left-hand plot shows the output of $\lnk$ for each input radius; the effect of the polynomial interpolant is limited to a narrow window near the upper and lower constraints. The plot on the right shows several choices of $\lnk$ applied to a triangular inner boundary. Even though the triangle is well within the minimum and maximum radius constraints, the standard functions nevertheless change the boundary quite a bit, while the polynomial interpolant lies directly on top of the triangle. 

\subsubsection{Mesh Generation and Finite Element Solver} \label{sec:fem}

To mesh the domain, one point for each B-spline used in \cref{sec:bspline} was added in {\tt Gmsh} to represent the inner and outer boundaries. A Bspline was then drawn through each set of points, a curve loop was added for each, and then a plane surface was created to represent the area between them. We then generated a mesh of the domain with quadratic elements using a target mesh size (\texttt{lc} parameter) of $0.03$. Integration of {\tt Gmsh} into the Julia code was conducted using the \texttt{Gmsh} Julia package (\url{https://juliapackages.com/p/gmsh}).

When subdomains/sensors were required as in \eqref{eq:Vol:ave:msr} or \cref{sec:ex:1}, the subdomains were generated in a way that shares the same points and edges along the subdomain boundaries. Example meshes with subdomains are shown in \cref{fig:vort:samples}.

\begin{figure}[ht!]
	\centering
	\includegraphics[width=0.45\linewidth]{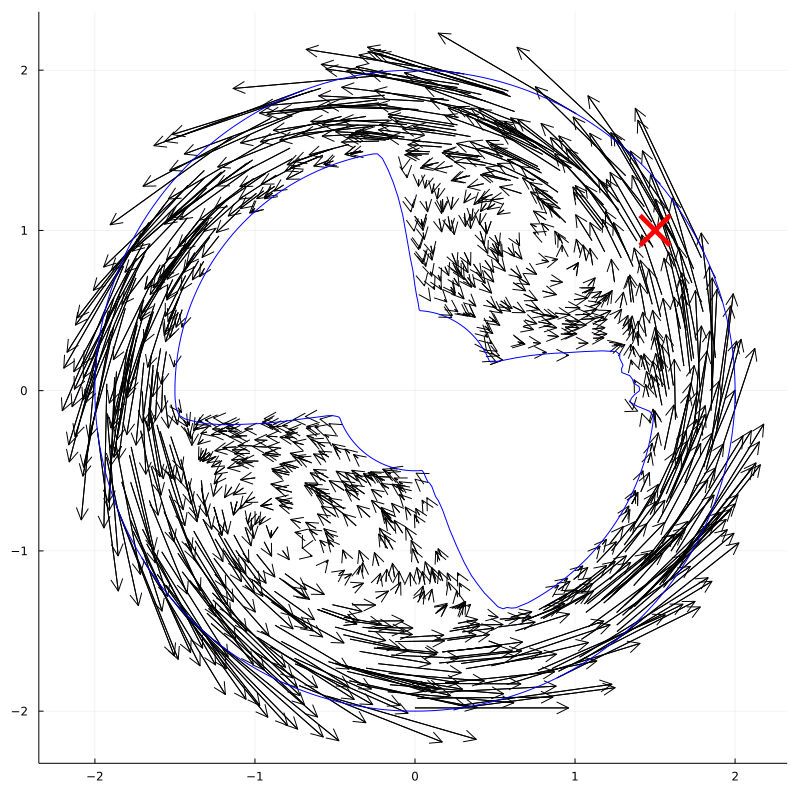}
	\hfill
	\includegraphics[width=0.51\linewidth]{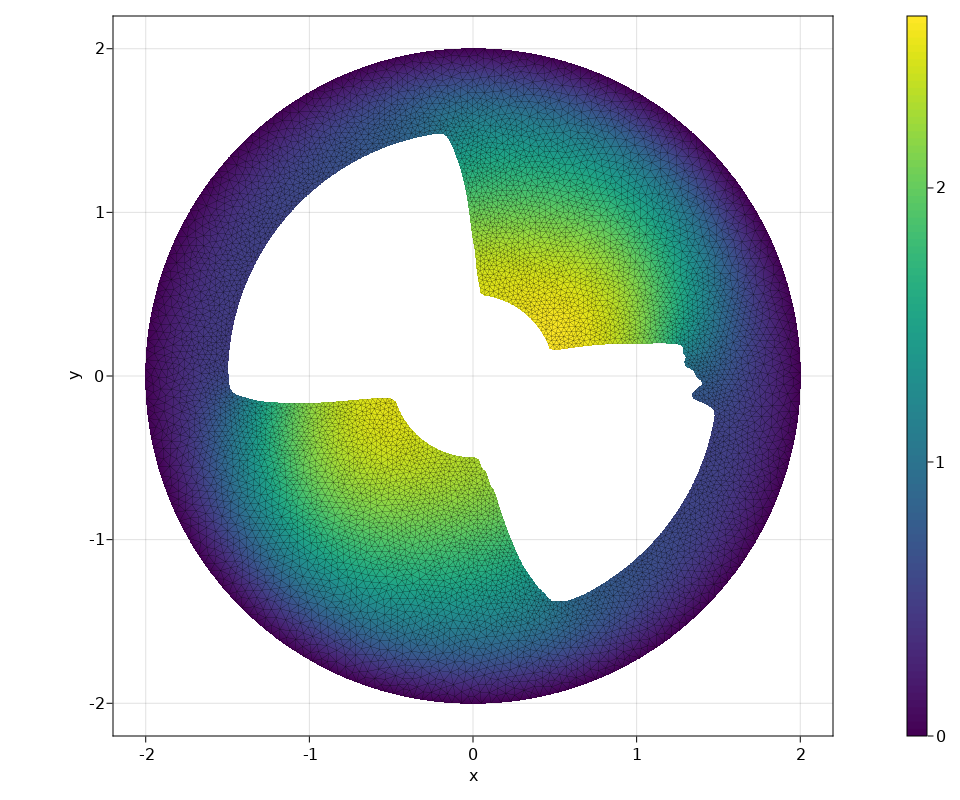}
	\caption{Solutions of the coupled Stokes and
          advection-diffusion PDEs. Left: Quiver plot of the solution
          to the Stokes PDE. The red X marks the point at which the
          scalar is injected into the system. Right: Contour plot of
          the solution to the advection-diffusion equation associated
          with the Stokes flow in the left-hand
          plot.}\label{fig:solutions}
\end{figure}

The triangulation provided by {\tt Gmsh} defines piecewise quadratic approximations to the velocity and advected scalar fields as well as a piecewise linear pressure field.  The steady Stokes equations \eqref{eq:Sk:bulk}, \eqref{eq:Sk:bnd} are thus approximated using a mixed formulation with the Taylor-Hood finite element pair, cf.~\cite{gunzburger1989FiniteElementMethods}.  A penalty method is used to regularize the equations, replacing the constraint $\nabla\cdot\bU=0$ with $\nabla\cdot\bU+\varepsilon p=0$.  Approximating the velocity $\bU$ with quadratic elements of nominal size $h$ and using the penalty parameter $\varepsilon$, it is known (cf.~\cite{gunzburger1989FiniteElementMethods}) that the Galerkin finite element approximation, $\bU_\varepsilon^h$ satisfies the estimate (in the usual Sobolev norm $\mathcal{H}^1(\DD_b)$)
\begin{displaymath}
  \| \bU - \bU_\varepsilon^h \|_{{\cal H}^1} \leq C\left( h^2 + \varepsilon\right).
\end{displaymath}
Accordingly, we chose a penalty parameter of $\varepsilon=0.001$ to have a similar contribution to the velocity error as the discretization error since $h\approx0.03$ for our choice of element sizes.  This was implemented using software adapted from routines found in \url{https://github.com/jborggaard/FEMfunctions.jl}.
\cref{fig:solutions} shows plots of the solution to the coupled Stokes and advection-diffusion equations for an example inner boundary.

\section{Example Problems}
\label{sec:example:problems}

This section provides numerical examples where we compute the approximate posterior measure for each of the three different observation types described in \cref{sec:obs:procedure}. In the first, we consider only the Stokes equation, seeking inner boundaries that yield flows matching given vorticity values in certain subregions as in \eqref{eq:Vol:ave:msr}. In the second and third examples, we couple the Stokes and advection-diffusion equations and seek shapes that yield certain features in the scalar field $\pS$ as in \eqref{eq:scalar:var:msr} and \eqref{eq:scalar:sect:var:msr}. The problem parameters that are common across all three examples are summarized in \cref{tab:stokes:parameters}. All numerical results were obtained using the research computing clusters provided by Advanced Research Computing at Virginia Tech (\url{https://arc.vt.edu}).

\begin{table}[htbp]
  {\footnotesize
  \caption{Parameter choices for the Stokes Problem.} \label{tab:stokes:parameters}
  \centering
  \begin{tabular}{|>{\raggedright}p{0.25\textwidth}|p{0.3\textwidth}||>{\raggedright}p{0.25\textwidth}|p{0.075\textwidth}|} \hline
    Parameter & Value &
    Parameter & Value \\
    \hline\hline
    Prior, $\mpr$ \eqref{eq:post:form} & $\bdf_{2k-1},\bdf_{2k} \sim N(0,k^{-2s-1})$ &
    Range for $\lnk$, $\epsilon$ \eqref{eq:clamp:poly} & $0.1$  \\\hline
    Radius Constraints \eqref{eq:com:domain} & $\ir=0.5,\,\mr=1.5,\,\oR=2$ &
    Mean radius, $\bdf_0$ \eqref{eq:b:dom} & $1.0$ \\\hline
    Sampling dimension, $K$ & 320 &
    Number of B-splines \eqref{eq:bSpl:basis} & 160 \\\hline
    Angular Velocity, $\rtrt$ \eqref{eq:Sk:bnd} & 10 &
    Diffusion, $\kappa$ \eqref{eq:AD:steady} & $1$ \\\hline
    Source, $\src(\x)$ \eqref{eq:AD:steady} & $4 \exp[-(\x-\x_0)^2/100]$, $\x_0=(1.5,1)$ &
    Kinematic viscosity, $\visc$ \eqref{eq:Sk:bulk} & $0.001$
    \\\hline
  \end{tabular}
  }
\end{table}

\subsection{Example 1: Targeting Vorticity}
\label{sec:ex:1}
In this example, we try to generate vorticity of certain magnitude in subregions of the domain; we call these locations the ``sensors'' below. Eight circular sensors are placed midway between the maximum inner boundary $\mr$ and the outer boundary $\oR$ and evenly spaced radially (every $45\degree$) around the domain. The observations are then computed as in \eqref{eq:Vol:ave:msr}. These problem parameters specific to this example are summarized in
\cref{tab:stokes:param:1}.
\begin{table}[htbp]
  \centering
  {\footnotesize
  \caption{Parameter choices for Example 1.} \label{tab:stokes:param:1}
  \begin{tabular}{|>{\raggedright}p{0.35\textwidth}|p{0.55\textwidth}|} \hline
    Parameter & Value \\
    \hline\hline
    Observations, $\Obs$ (\cref{sec:obs:procedure}) & Vorticity \eqref{eq:Vol:ave:msr} at $8$ sensors (\cref{fig:vort:samples}) \\\hline
    Data, $\Data = (\data_1, \dots, \data_8)$ \eqref{eq:stat:inverse:model} & $\data_1=\data_7=30;~\data_3=\data_8=50; \text{ and } \data_j=40 \text{ otherwise}$ (\cref{fig:vort:quantiles}, Left)\\\hline
    Noise, $\noise$ \eqref{eq:stat:inverse:model} & $N(0,\sigma_{\noise}^2 I)$, $\sigma_\noise=1.0$ \\\hline
    Sensor locations, $\x_j$ \eqref{eq:Vol:ave:msr} & $\x_j= \left( 1.75\cos \pth_j, 1.75\sin \pth_j \right), \pth_j=\pi j / 4 \text{ for }j=1,\ldots,8$ \\\hline
    Sensor radii, $r_j$ \eqref{eq:Vol:ave:msr} & $r_1=\dots=r_8=0.1$ \\\hline
    Sobolev regularity $s$ of the boundary  \eqref{eq:sob:sp} & $s =1.25$ \\\hline
  \end{tabular}
  }
\end{table}

\begin{figure}
    \centering
    \includegraphics[width=0.7\textwidth]{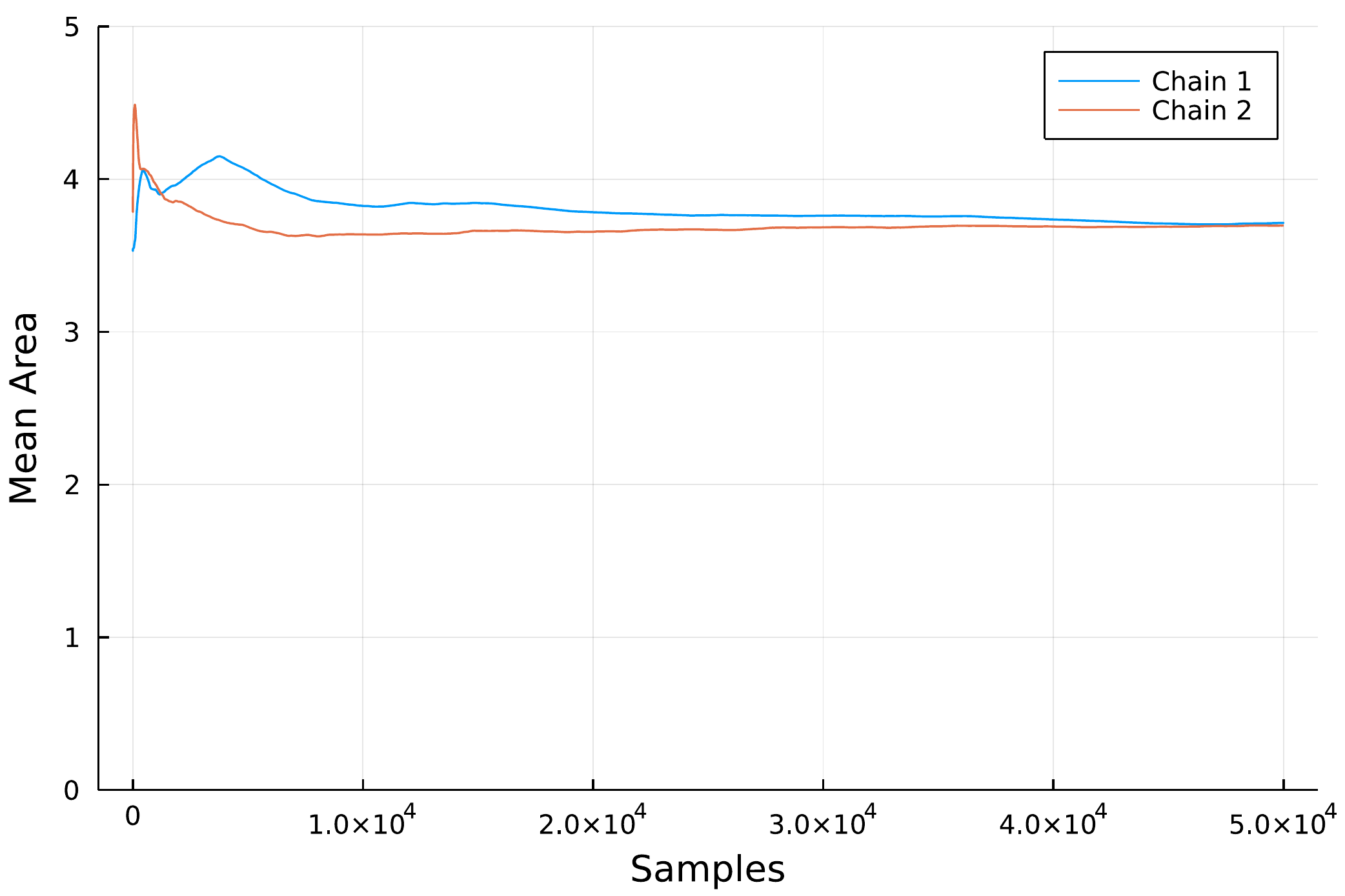}

    \caption{Running average of the area enclosed by the inner boundary, by chain.}
    \label{fig:vort:converge}
\end{figure}

To compute the approximate posterior measure, we ran two pCN chains with $50,000$ samples (approximately two days' worth of run time) each. The final value of $\rho$ (see \cref{alg:pcn}) was $0.9993$ and the resulting acceptance rate was $25.9\%$. To confirm that the chains' results were similar, we computed the $\hat{R}$ diagnostic from \cite[Section 11.4]{gelman2014bayesian}, a common metric for measuring convergence of MCMC chains. Using the area enclosed by the inner boundary as our estimand, we found $\hat{R}=1.002$, indicating convergence. \cref{fig:vort:converge} shows a running average of the area of the inner boundary as the chains progress.

\begin{figure}
    \centering
    \includegraphics[width=0.44\textwidth]{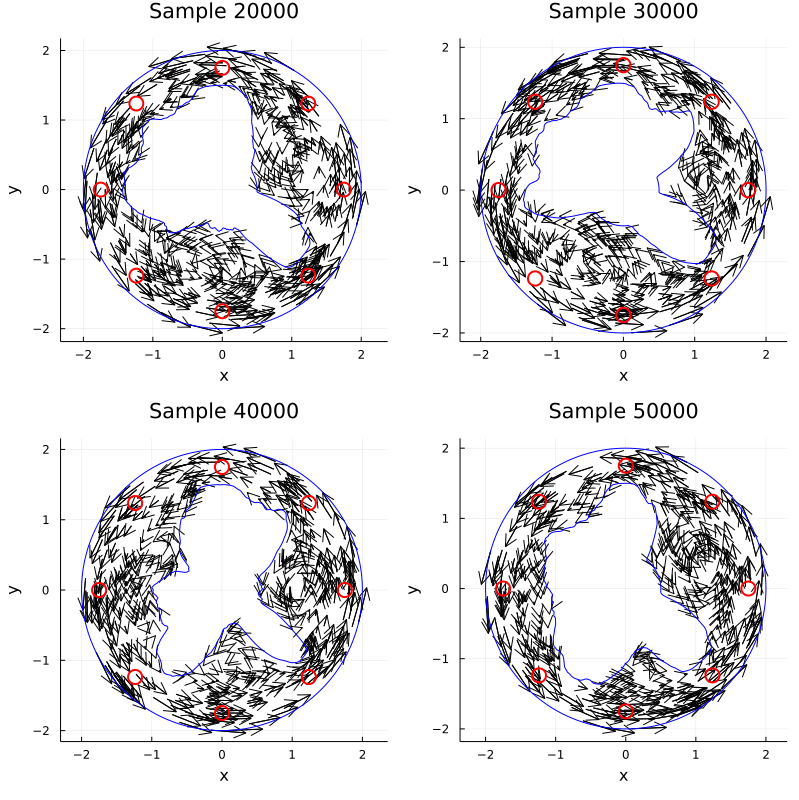}
    \hfill
    \includegraphics[width=0.54\textwidth]{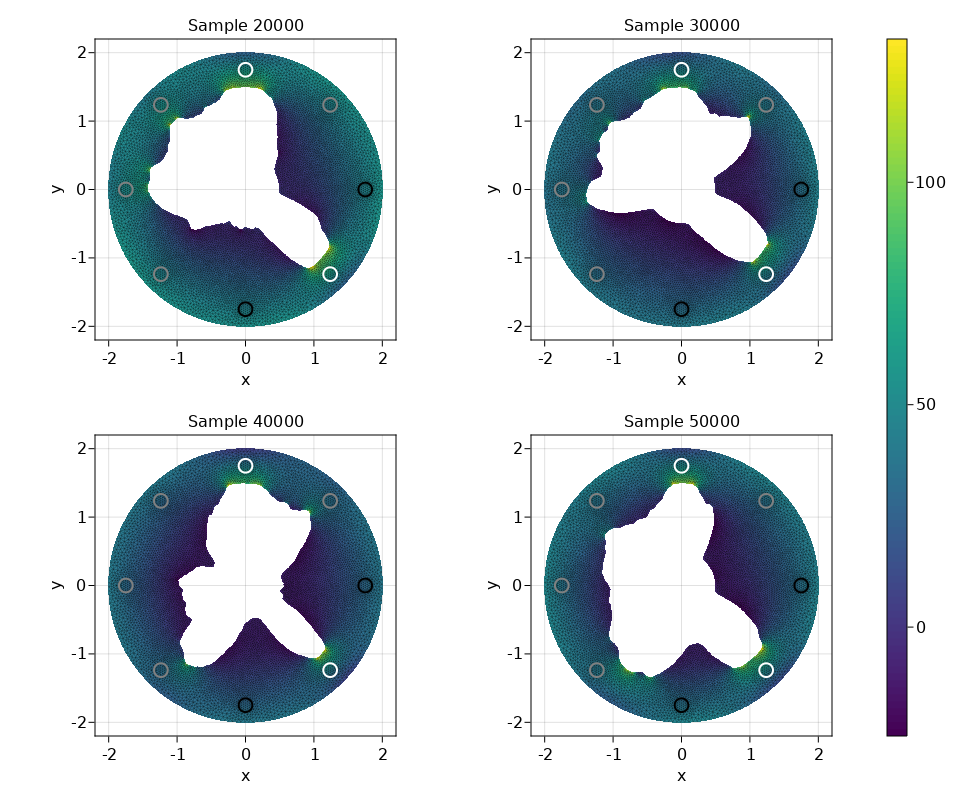}
    \caption{Example posterior MCMC samples for the vorticity problem. Left: Quiver plot of the resulting flow (sensors in red). Right: Heatmap of vorticity; the target data is high vorticity at the top and bottom-right sensors (white circle in the plot), low vorticity at the right and bottom sensors (black), and medium vorticity at other sensors (gray).  See \cref{tab:stokes:param:1} for the exact values of the this target data.}    \label{fig:vort:samples}
\end{figure}

\cref{fig:vort:samples} shows four example samples taken from one of the chains. Each of these samples does a good job of matching both the data $\Data$ and the structure imposed by the prior $\mpr$. The left plot shows quiver plots of the associated Stokes flows, while the right plot shows heatmaps of flow vorticity. In the latter plot, we have colored the sensors by target data value, i.e., white/gray/black indicate high/medium/low target vorticity, respectively. We see that the shapes tend to have the radius of the interior boundary be largest near the top and bottom-right sensors where the vorticity data takes its highest values. Such data points seem to push the fixed inner boundary close to the moving outer boundary to induce more vorticity in the flow. Similarly, the radius is smallest -- close to the minimum value $\ir$ -- at the bottom and right sensors where the data dictates that the vorticity should be smaller.

\begin{figure}
    \centering
    \includegraphics[width=0.40\textwidth]{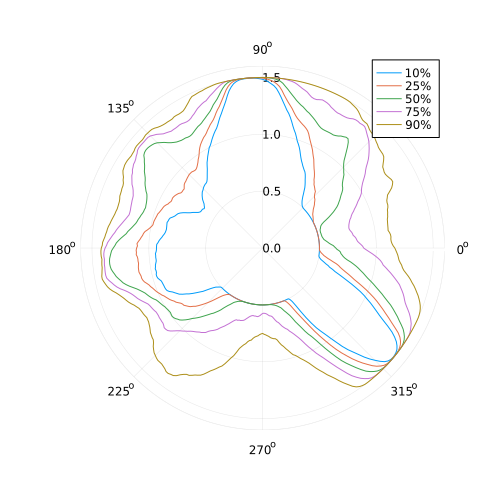}
    \hfill    
    \includegraphics[width=0.58\textwidth]{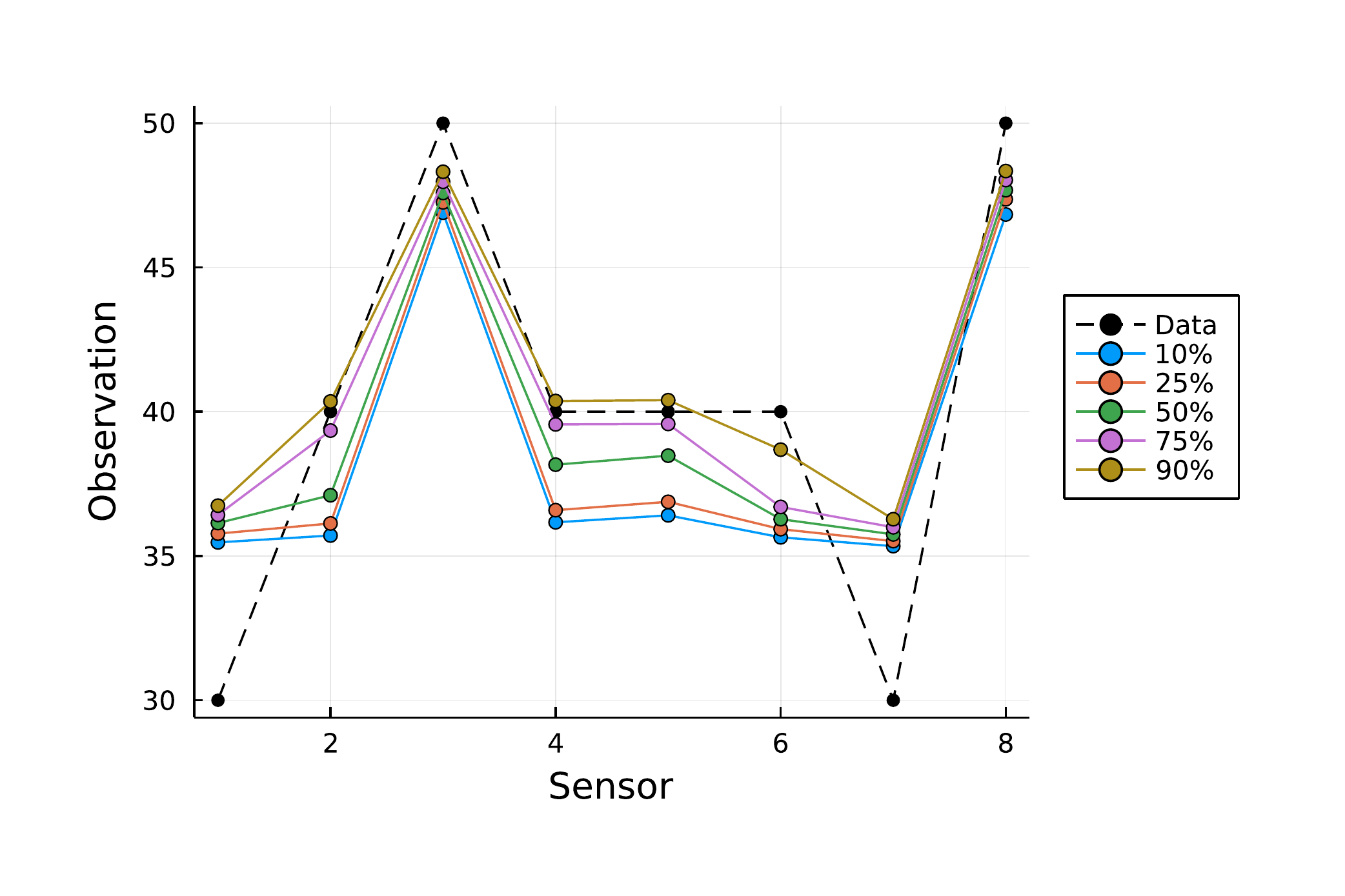}
    \caption{Quantiles for the vorticity problem. 
    Left: Radius quantiles.
    Right: Observation quantiles and target data. Sensors are numbered counterclockwise starting from the the right-hand side of the region. 
    }
    \label{fig:vort:quantiles}
\end{figure}

While \cref{fig:vort:samples} shows a few solutions, the power of the Bayesian approach is that the ``answer'' is a probability distribution indicating some degrees of freedom in the design. We illustrate some features of this distribution in \cref{fig:vort:quantiles}, which shows quantiles of the posterior. The left-hand plot shows the quantiles of the radii. To generate this figure, for each angle from $0\degree$ to $360\degree$, we compute the radius for each sample, and then compute the quantiles for each angle; we then plot each quantile on the polar plot. The quantiles show that virtually all designs maximized the inner radius at the top and bottom right of the domain where the vorticity was specified to be highest; similarly, most designs minimized the inner radius at the bottom and right of the domain where vorticity was to be lowest. The design appears to have more degrees of freedom near the other four sensors (the left half of the domain and at top right); in these regions the vorticity was set to a middle value and the radii seem to vary more. The right-hand plot of \cref{fig:vort:quantiles} shows the quantiles of the resulting observations \eqref{eq:Vol:ave:msr} plotted against the data $\Data$, indicating which values the designs could achieve and which they could not. Here we see that the designs achieved the overall shape of the data, with peaks that match the high values and some variation near the middle vorticity values. The ``low vorticity'' target of $30$ was perhaps overly aggressive as the designs could not drop vorticity below roughly 35; 
however the low target values and constant observational noise deviation $\sigma_{\noise}$ forced the observations at these locations to cluster close to the lower bound. Overall, these figures illustrate which features of the shape are dictated by the optimization problem and which are more free to vary.

\subsection{Example 2: Trapping the Scalar}\label{sec:svglobal}

In this example, we couple the Stokes flow  \eqref{eq:Sk:bulk}--\eqref{eq:Sk:bnd} to the (steady) advection-diffusion equation  \eqref{eq:AD:steady}--\eqref{eq:AD:steady:BC} as described above, and seek to identify shapes that will achieve certain behavior in the scalar $\pS$. As noted in \cref{sec:obs:procedure}, a key measure of mixing is the scalar variance \eqref{eq:scalar:var:msr}; in this case we target global average scalar variance of $0.8$, which is a relatively high value given the other problem parameters, indicating that the scalar has been ``trapped'' in a subregion of the domain. This choice of problem parameters is summarized in \cref{tab:stokes:param:2}. Several MCMC chains were run for 25,000 samples each and the results compared to ensure that they were consistent with one another.

\begin{table}[htbp]
  \centering
  {\footnotesize
  \caption{Parameter choices for Example 2.} \label{tab:stokes:param:2}
  \begin{tabular}{|>{\raggedright}p{0.35\textwidth}|p{0.35\textwidth}|} \hline
    Parameter & Value \\
    \hline\hline
    Observations, $\Obs$ (\cref{sec:obs:procedure}) & Average scalar variance \eqref{eq:scalar:var:msr} \\\hline
    Data, $\Data$ \eqref{eq:stat:inverse:model} & $0.8$ \\\hline
    Noise, $\noise$ \eqref{eq:stat:inverse:model} & $N(0,\sigma_{\noise}^2 I)$, $\sigma_\noise=0.05$ \\\hline
    Sobolev regularity $s$ of the boundary \eqref{eq:sob:sp} & $1.0$ \\\hline
  \end{tabular}
  }
\end{table}

\begin{figure}[h]
    \centering
    \includegraphics[width=0.45\textwidth]{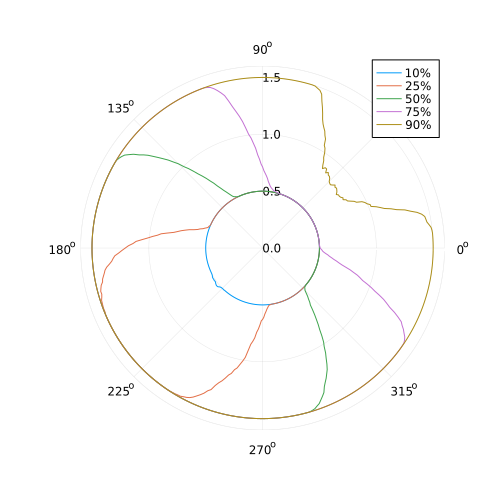}
    \hfill
    \includegraphics[width=0.51\textwidth]{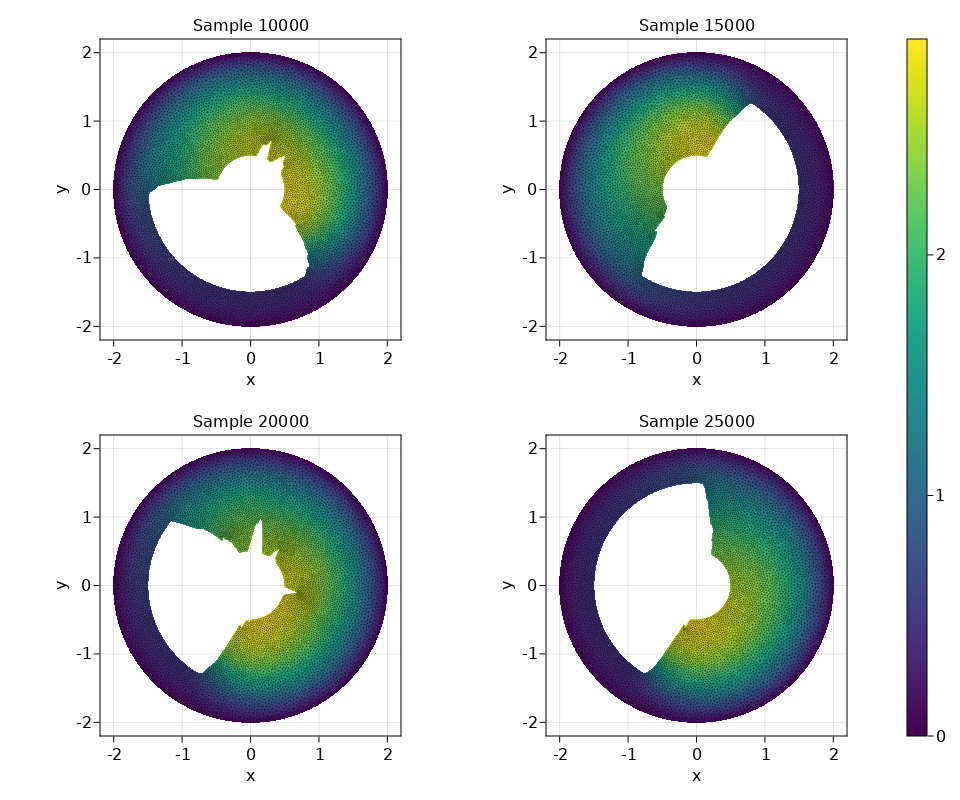}
    \caption{Left: Radius quantiles for the global scalar problem. Right: Plot of scalar $\pth$ for four posterior MCMC samples for the global scalar problem.}
    \label{fig:svglobal:radii:quantiles}
\end{figure}

The left-hand plot of \cref{fig:svglobal:radii:quantiles} shows quantiles of the shape radii by angle. We see that the radius is nearly always the minimum at $45\degree$ and the maximum at $225\degree$, but other angles show some variation.
 The right-hand plot of \cref{fig:svglobal:radii:quantiles} shows the boundary shape and the scalar field for four samples drawn from the chain; these samples imply that the posterior consists of ``bell'' shapes with some degree of freedom in their orientation.
 
\begin{figure}
    \centering
    \includegraphics[width=0.7\textwidth]{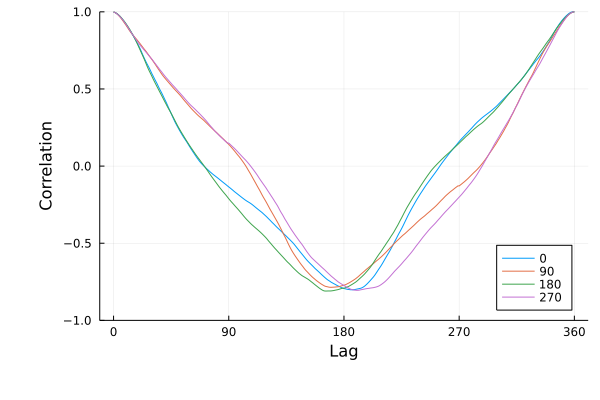}
    \caption{Radius correlations by angle and lag for the global scalar problem.}
    \label{fig:svglobal:radii:corr}
\end{figure}

To understand this rotation and the variance in the radius quantiles, we can compute the correlations between radii at different angles; \cref{fig:svglobal:radii:corr} shows the correlation of the radius at $0$, $90$, $180$, and $360$ degrees with other radii at various lags (offsets). The plot shows that radii are strongly correlated -- radii tend to go down when at the opposite side of the domain go up, and vice versa. This anticorrelation shows that the bell shape is rotating, which reflects a degree of freedom in the problem -- because scalar variance is a rotation-independent metric, the inverse problem has a fair amount of latitude in how to orient the bell. The only rotation-dependence, and therefore the only primary constraint to the orientation, is the point where the scalar is injected to the system; if the ``thick'' part of the bell coincides with the injection point, the scalar is pushed near the boundary where it is removed from the system, driving down scalar variance. This explains why the posterior does not appear to include shapes with large radius near the injection point. 

One benefit of the statistical approach to the shape estimation problem is that the posterior distribution makes clear in this case that there is a degree of freedom in choosing a shape that approximates the target data. In the next example, we define observations that are rotation dependent and therefore place greater constraint on the shape orientation.

\subsection{Example 3: Adding Rotation Dependence}\label{sec:svsector}
For our final example, we adjust the previous observational setup to introduce a rotation-dependence, thereby better constraining the desired shape. We do this by replacing average scalar variance across the domain \eqref{eq:scalar:var:msr} as our observation operator with scalar variance by quadrant. That is, we use sectoral scalar variance \eqref{eq:scalar:sect:var:msr} where the sectors $\qS_j(\DD_b), j=1,\dots,4$ as in \eqref{eq:sectors} are defined with angles $\pth_j = \pi j / 2, j=0,\dots,4$. We require that scalar variance be high (target 0.4) in the top right and bottom left quadrants and low (target 0.0) in the top left and bottom right quadrants. The parameters in the inference are otherwise the same as in the previous examples. These parameter choices are summarized in \cref{tab:stokes:param:3}. 
\begin{table}[htbp]
  \centering
  {\footnotesize
  \caption{Parameter choices for Example 3.} \label{tab:stokes:param:3}
  \begin{tabular}{|>{\raggedright}p{0.35\textwidth}|p{0.25\textwidth}|} \hline
    Parameter & Value \\
    \hline\hline
    Observations, $\Obs$ (\cref{sec:obs:procedure}) & Sectoral scalar variance \eqref{eq:scalar:sect:var:msr} \\\hline
    Data, $\Data = (\data_1, \data_2, \data_3, \data_4)$ \eqref{eq:stat:inverse:model} & $\data_1=\data_3=0.4,\, \data_2=\data_4=0$ \\\hline
    Noise, $\noise$ \eqref{eq:stat:inverse:model} & $N(0,\sigma_{\noise}^2 I)$, $\sigma_\noise=0.05$ \\\hline
    Sobolev regularity  $s$ of the boundary \eqref{eq:sob:sp} & $1.0$ \\\hline
  \end{tabular}
  }
\end{table}

\begin{figure}
    \centering
    \includegraphics[width=0.45\textwidth]{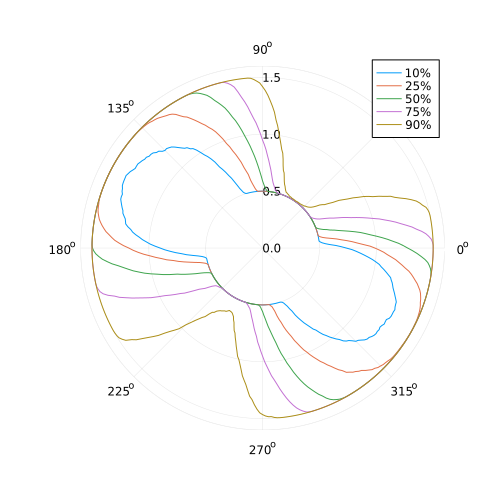}
    \hfill
    \includegraphics[width=0.51\textwidth]{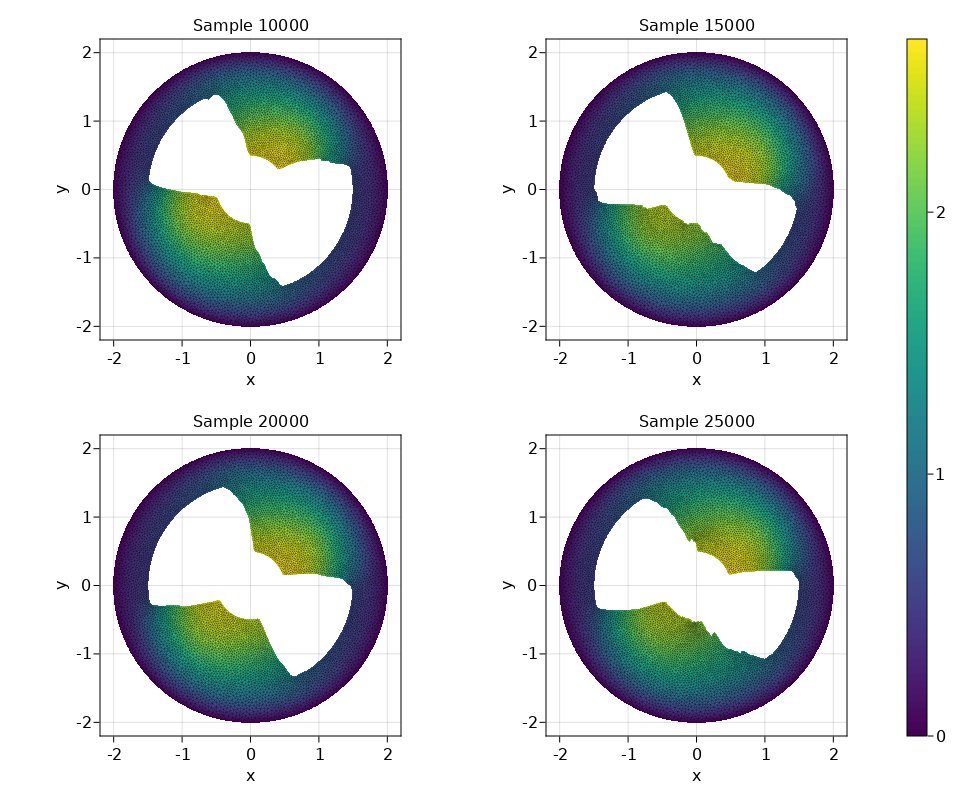}
    \caption{Left: Radius quantiles for the sectoral scalar problem. Right: Plot of scalar $\pth$ for four  MCMC samples for the sectoral scalar problem.}
    \label{fig:svsector:radii:mle}
\end{figure}

The left hand plot of \cref{fig:svsector:radii:mle} shows the radii quantiles by angle for the sectoral scalar variance problem. As with \cref{fig:svglobal:radii:quantiles}, we see clear constraints on the shape; in this case the rotor takes on a ``bowtie'' shape where the radius is consistently near the minimum at roughly $60\degree$ and $240\degree$ and near the maximum at approximately $150\degree$ and $330\degree$. 
We see this shape for the samples in the right hand plot of \cref{fig:svsector:radii:mle}, which also shows that the shape traps the scalar in two regions that align well with top right and bottom left quadrants. 
By contrast, the wider lobes of the bowtie serve to minimize the amount of scalar in the top left and bottom right quadrants. The lobes appear to have small degrees of freedom to move and stretch, but seem to be much more constrained than the shapes in \cref{sec:svglobal}.

\begin{figure}
    \centering
    \includegraphics[width=0.7\textwidth]{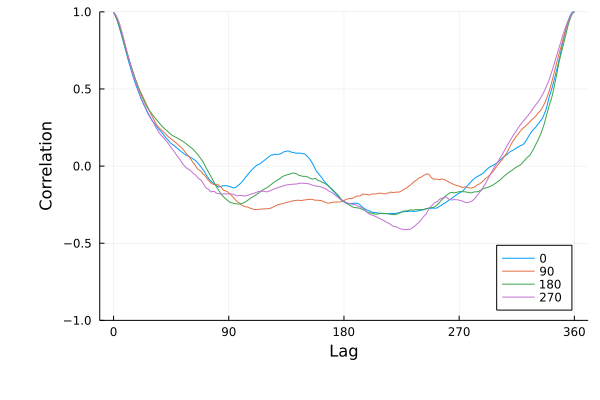}
    \caption{Radius correlations by angle and lag for the sectoral scalar problem.}
    \label{fig:svsector:radii:corr}
\end{figure}

In the radii correlations \cref{fig:svsector:radii:corr}, however, we see a clear contrast with Example 2 (\cref{fig:svglobal:radii:corr}). The plot shows minimal correlations between the radii at angles far apart from each other. This shows that while the radii do have some freedom to vary, as seen at the edges of the lobes in \cref{fig:svsector:radii:corr}, this variation is less a result of a rotational freedom and more because of other ranges of motion -- e.g., a given lobe growing or shrinking.

\section{Outlook} \label{sec:conclusions}

A number of interesting avenues for future work immediately present themselves which build
on and compliment the methodology we developed here.  Firstly, while the numerical 
result in \cref{sec:example:problems} strongly suggests that the forward map 
defining the posterior measure in  \eqref{eq:post:form} has desirable analytical
properties, a complete and systematic analysis of $\mathcal{G}$ remains for
future research endeavors.  Here a complete picture would include 
continuous dependence of the posterior on observed data as advocated in 
\cite{dashti2017bayesian}; see also \cite{cohen2018ShapeHolomorphyStationary} for a related result deriving desirable properties on the map from boundary parameters to PDE solution for Stokes and Navier-Stokes flows.  For such an effort we anticipate that existing machinery from
e.g.~\cite{delfour2011shapes} may provide core elements of the required analysis 
but it is by no means obvious at the outset that this represents a routine exercise.
Another significant and related research task is to derive adjoint methods
for the computation of gradients in these problems.  Such result would open
up the possibility of employing HMC and MALA type sampling methods for the resolution of
posteriors.  One could then consider the need for consistent adjoints~\cite{Hartmann2007AdjointConsistency,hiptmair2015ComparisonApproximateShape} and their impact on convergence of the posteriors. 

As we alluded to above, particularly in \cref{sec:mcmc}, the example problems 
presented in this work provide a series of challenging benchmark test cases for 
the development of and the optimal tuning of algorithmic parameters in infinite dimensional 
sampling methods.    While the computations herein and in a number of complementary recent works
clearly illustrate the power and potential of methods developed in e.g. 
\cite{Beskosetal2011,cotter2013mcmc, glatt-holtz2020acceptreject, glattholtz2022on} much more 
experimental work is warranted to better understand the relative efficacy of the different possible approaches
at hand.
In particular the estimation problems considered here provide a setting where the 
computation of gradients is a particularly onerous task.
Here our recent work in \cite{glatt-holtz2020acceptreject} provides a systematic
approach to using \emph{approximations} of the gradient of the forward map.  
Meanwhile in \cite{glattholtz2022on} we developed a novel gradient free method mpCN generalizing the 
pCN algorithm which produces a cloud of proposals at each algorithmic step.  A series of preliminary
case studies of this new mpCN algorithm were carried out in  \cite{glattholtz2022on} one of which used the example in
\cref{sec:svsector} as a test case.

Finally let us mention that a variety of other boundary shape estimation problem present themselves
for which the numerical machinery and statistical framework we have build here
can be employed with modest adjustments to our code base. These could include more complicated fluids cases more relevant to  applications areas, such as for example by adding a boundary drag term to account for the torque required to rotate the outer boundary. 
Meanwhile, one classical PDE inverse problem 
with seemingly no precedent in the Bayesian literature is the estimation
a two dimensional domain shape from the associated eigenvalue structure 
of the Dirichlet-Poisson eigenvalue problem. In effect we might pose the question: 
does a Bayesian set-up  prove useful for regularizing the ill-posed problem of 
`hearing the shape of a drum.'  In each of these cases, adopting a statistical approach as we describe here might illuminate new features of the shape estimation problem.

\section*{Acknowledgements}

Our efforts are partially supported by the National Science Foundation under grants NSF-DMS-1819110, NSF-CMMI-2130695 (JTB), NSF-DMS-1816551, NSF-DMS-2108790 (NEGH), 
 NSF-DMS-2108791 (JAK). The authors acknowledge Advanced Research Computing at Virginia Tech (\url{https://arc.vt.edu}) for providing computational resources and technical support that have contributed to the results reported within this paper.

\addcontentsline{toc}{section}{References}
\begin{footnotesize}
\bibliographystyle{alpha}
\bibliography{bib}
\end{footnotesize}

\newpage
\begin{multicols}{2}

\noindent
Jeff Borggaard\\
{\footnotesize
Department of Mathematics\\
Virginia Tech\\
Web: \url{https://personal.math.vt.edu//borggajt/research/}\\
Email: \url{jborggaard@vt.edu}}\\ [.2cm]

\noindent
Nathan E. Glatt-Holtz\\ {\footnotesize
Department of Mathematics\\
Tulane University\\
Web: \url{http://www.math.tulane.edu/~negh/}\\
Email: \url{negh@tulane.edu}} \\

\columnbreak

\noindent 
Justin Krometis\\
{\footnotesize
National Security Institute\\
Virginia Tech\\
Web: \url{https://krometis.github.io/}\\
Email: \href{mailto:jkrometis@vt.edu}{\nolinkurl{jkrometis@vt.edu}}} \\[.2cm]

\end{multicols}

\end{document}